\documentclass[reqno,12pt,a4paper]{amsart}
\usepackage{amsmath,amssymb,amsthm,amscd,amsfonts,hyperref,amstext}
\usepackage[english]{babel}
\usepackage[mathscr]{eucal}

\numberwithin{equation}{section}

\renewcommand{\a}{\alpha}
\renewcommand{\b}{\beta}

\renewcommand{\d}{\delta}
\newcommand{\D}{\Delta}
\newcommand{\e}{\epsilon}

\newcommand{\z}{\zeta}

\newcommand{\m}{\mu}
\newcommand{\n}{\nu}
\newcommand{\x}{\xi}

\newcommand{\s}{\sigma}

\renewcommand{\t}{\tau}
\newcommand{\f}{\phi}

\newcommand{\F}{\Phi}

\newcommand{\p}{\psi}

\renewcommand{\O}{\Omega}

\newcommand{\C}{{\mathbb C}}
\newcommand{\R}{{\mathbb R}}
\newcommand{\Z}{{\mathbb Z}}

\newcommand{\Ab}{{\mathbf A}}

\newcommand{\Hb}{{\mathbf H}}

\newcommand{\Ac}{{\mathcal A}}
\newcommand{\Bc}{{\mathcal B}}

\newcommand{\Ec}{{\mathcal E}}
\newcommand{\Fc}{{\mathcal F}}

\newcommand{\Hc}{{\mathcal H}}

\newcommand{\Lc}{{\mathcal L}}

\newcommand{\Sc}{{\mathcal S}}

\newcommand{\supp}{\hbox{{\rm supp}}\,}

\DeclareMathOperator{\im}{{\rm Im}\,}
\DeclareMathOperator{\re}{{\rm Re}\,}
\DeclareMathOperator{\rank}{rank}

\newcommand{\card}{\operatorname{card}}

\newcommand{\psupp}{\operatorname{psupp}}
\newcommand{\singsupp}{\operatorname{sing supp}}

\newtheorem{theorem}{Theorem}[section]
\newtheorem{proposition}[theorem]{Proposition}
\newtheorem{lemma}[theorem]{Lemma}
\newtheorem{corollary}[theorem]{Corollary}

\theoremstyle{definition}

\theoremstyle{remark}

\newtheorem{example}[theorem]{Example}

\begin{document}

\title[Finite rank Toeplitz operators ]{Finite Rank Toeplitz Operators: Some Extensions of  D.Luecking's Theorem }
\author[Alexandrov]{Alexey Alexandrov}
\address[A. Alexandrov]{Petersburg Department of Steklov Institute of Mathematics,
 Russian Academy of Sciences. 27, Fontanka, St. Petersburg, 191023,  Russia}
\email{alex@pdmi.ras.ru}
\author[Rozenblum]{Grigori Rozenblum}
\address[G. Rozenblum]{Department of Mathematics \\
                          Chalmers University of Technology \\
                          and Department of Mathematics University of Gothenburg \\
                          Chalmers Tv\"argatan, 3, S-412 96
                           Gothenburg
                          Sweden}
\email{grigori@math.chalmers.se}

\begin{abstract}
The recent theorem by D.Luecking about finite rank Bergman-Toeplitz operators is extended to weights
being distributions with compact support and to the spaces of harmonic functions.
\end{abstract}
\keywords{Bergman spaces, Bargmann spaces, Toeplitz
operators}
\date{}

\maketitle


\section{Introduction and the main result}
\label{intro} Toeplitz operators play an important role in many branches of analysis. A significant
recent development in the theory of such operators is related to the proof, given by D.Luecking
\cite{Lue2}, of the finite rank conjecture. Let $\Bc^2$ be the Bergman space of $L_2$-functions
analytical in a domain $\O\in\C^1$ and  $P$ be the orthogonal  projection in $L_2(\O)$ onto $\Bc_2$.
For a regular complex Borel measure $\mu$ with compact support, the Toeplitz operator with weight
$\mu$,
\begin{equation}\label{1:Toepldef}
    u\mapsto T_\mu u=Pu\mu , u\in \Bc^2,
\end{equation}
can be correctly defined. According to the finite rank conjecture, if $T_\mu$ has finite rank then
the measure is a finite combination of point masses, exactly as many as the rank is.  The nontrivial
past of this conjecture is described in \cite{Lue2}, \cite{RShir}. Immediately after the preprint
containing the proof appeared, an activity developed in extending and applying  this result. On the
one hand, the theorem by Luecking was extended to the multi-dimensional case, see \cite{BRChoe1},
\cite{RShir} (by different methods). On the other hand, interesting application to the theory of
Toeplitz operators appeared, see \cite{BRChoe2}, \cite{Cuc}, \cite{Le} \cite{Le1}, as well as in
Function Theory, see  \cite{BRChoe1}.  The finite rank result turns out to be useful also in
Mathematical Physics, more exactly, to the spectral analysis of the perturbed Landau Hamiltonian, see
\cite{RaiWar}, as well as the discussion and further references in \cite{RShir}.

A number of natural questions arise around Luecking's
theorem. First, it is interesting to find out whether
the finite rank property still holds when the
analytical Bergman space is replaced by some other,
also closed in $L_2$,  space of smooth functions. In
\cite{BRChoe1} such a generalization was found for the
space of $n$--harmonic functions in a domain in
$\C^n$, and in \cite{Le} the finite rank property was,
in the complex dimension 1, extended to the
$L_2$--closed span of  certain, not too sparse, sets
of monomials $z^{n_k}$, $n_k\in \Z_+$. At the same
time, for the problems arising in Mathematical
Physics, it is important to generalize the results to
the case when the weight  measure $\mu$ is replaced by
a distribution with compact support.

In the present paper we deal with these questions. First, in the complex dimension 1, for the
analytical Bergman space,  we describe the procedure of reducing the finite rank problem for a
distribution to the same problem for an absolutely continuous measure $\mu$, which is already taken
care of. Thus, the finite rank problem finds its solutions also for distributional weights. We note
that the reduction above seems to be necessary. The initial proof with measure weight was critically
based upon a lemma on the density of symmetric  polynomials of a special form in the space of
symmetric continuous functions of many variables, proved by an ingenious use of the Stone-Weierstrass
theorem. The distributional case requires a similar density result in the space of differentiable
functions, where no proper analogy  of the Stone-Weierstrass theorem exists. Moreover, the density
result itself turns out to be wrong for differentiable functions. We present an example demonstrating
this. Therefore our approach seems to be at the moment the only one able to treat the distributional
case.

The results on finite rank problem for distributional weights are further extended to the
multi-dimensional case. We use a modification of the induction on dimension  presented in
\cite{RShir}. It seems that  the  approach to  proving the multi-dimensional Luecking's theorem,
proposed in \cite{BRChoe1}  using Stone-Weierstrass argument would not work for distributions, by the
reasons given above.

Finally, we consider the finite rank problem for the Bergman space of harmonic functions. The result
follows immediately from the one in the analytical case in an \emph{even} dimension, since the space
of harmonic functions contains the space of $n$-harmonic functions, where the finite rank property is
an obvious consequence  of the one in  the analytical case, see \cite{BRChoe1}. Quite different is
the situation in an \emph{odd} dimension ($\ge3$), where no direct coupling of harmonic functions to
analytical ones exists. Here we are able to handle only the case of a measure acting as weight, using
 a sort of dimension-reduction argument and some Harmonic Analysis technique. We give also an
example, not disproving the finite rank conjecture directly, but just hinting that the situation here  with
distributions might be considerably more delicate than the one with measures.

The results of the paper were obtained when the first author enjoyed the hospitality of the
Department of Mathematics of Chalmers University of Technology in Gothenburg, Sweden, supported by
the grant from the Swedish Royal  Academy of Sciences, for which he expresses his gratitude.

\section{Setting}\label{Sect2}
Let $\tau$ be a positive measure in a domain $\O\subset{\C}^d$ such that $0<\int_{{\C}^d}|P|\,d\tau<+\infty$
for every polynomial $P$ of the complex variables $(z_1,\dots,z_d)$, $P\not\equiv0$. We consider the
space $L^2(\O,\t)$ and the subspace $\Ac(\O,\tau)\subset L^2(\O,\t)$ consisting of analytical
functions. It is a closed subspace, and we denote by $P_{\Ac}(\O,\t)$ the orthogonal projection onto
$\Ac(\O,\tau)$. Further on, as soon as the domain and the measure are fixed, we suppress them in the
notations. The typical examples here are the Bergman spaces, for the case of a bounded $\O$ with
(say) Lebesgue measure, and the Fock-Bargmann spaces for $\O=\C^d$, $\t$ being the Gaussian measure.
The projection $P_{\Ac}$ is an integral operator with the reproducing kernel $P(z,w)$, infinitely
smooth, analytical in $z$ and anti-analytical in $w$ in the domain $\O$.

Let $F$ be a distribution with compact support in $\O$, $F\in \mathcal{E}'(\O)$. We denote by
$\langle F,\f\rangle$ the action of the distribution $F$ on the function $\f\in C^\infty(\O)$. Then,
for  $u\in\Ac(\O,\tau)$, the expression
\begin{equation}\label{distr.weight}
(T_Fu)(z)=\langle F,P(z,\cdot)u(\cdot)\rangle, \ z\in\O,
\end{equation}
defines an analytical function $(T_Fu)(z)\in\Ac(\O,\tau)$. The corresponding operator $u\mapsto T_Fu$
is a natural generalization of the Toeplitz operator $u\mapsto PFu$, $u\in \Ac(\O,\tau)$ for the case
when $F$ is a bounded measurable function with compact support in $\O$. The operator $T_F$ is bounded
in $\Ac$. Its sesquilinear form can be described as
\begin{equation}\label{QForm}
    (T_Fu,v)=\langle F,u\bar{v}\rangle, \ u,v\in\Ac.
\end{equation}

In the special case when the distribution $F$ is, in fact, a complex  Borel measure $\m$ with compact
support in $\O$, the operator $T_F$ can be described as
\begin{equation}\label{meas}
(T_Fu)(z)=\int_\O P(z,w)u(w)d\m(w),
\end{equation}
and the sesquilinear form is given by
\begin{equation}\label{QFormMeas}
    (T_Fu,v)=\int u\bar{v} d\m.
\end{equation}

Suppose that the operator $T_F$ has finite rank, $\rank(T_F)=m<\infty.$ This means, in particular,
that for any, finite or infinite, system of functions $f_\a\in \Ac$, the system of functions
$g_\a=T_Ff_\a$ is linearly dependent and $\rank\{g_\a\}\le m$. This is correct, in particular, if we
take as $f_\a$ the system of polynomials $f_\a=z^\a$, $\a=(\a_1,\dots,\a_d)\in(\Z_+)^d$ . Therefore
the infinite matrix
\begin{equation}\label{Matrix}
    \Ab_F=\left(a_{\a\b}\right),\ \ a_{\a\b}=(T_Fz^\a,z^\b)=\langle F,z^\a\bar{z}^\b\rangle
\end{equation}
has finite rank, $\rank(\Ab_F)\le m$. It is important that the matrix $\Ab_F$ does not depend on the
domain $\O$ or the measure $\t$, but it depends only on the distribution $F$. Of course, the rank of
$\Ab_F$ does not change if we make a unitary transformation of $\C^d$ with corresponding change of
complex coordinates.

 We notice also, following \cite{RShir}, that if $g$ is  function analytical and bounded
in some polydisk neighborhood of $\supp F$  and $F_{g}$ is the distribution $|g|^2F$ then
$\rank\Ab_{F_g} \le \rank \Ab_F$. To show this, we consider first  a polynomial $g_l$ of degree $l$.
The matrix $\Ab_{F_{g_l}}$ is obtained by building linear combinations of rows and columns of
$\Ab_F$, therefore the rank does not increase, $\rank\Ab_{F_{g_l}} \le \rank \Ab_F$. We pass to a
general analytical function $g$ using approximations by Taylor polynomials, convergent, together with
all derivatives, uniformly on any compact in the polydisk.

In a similar way, we consider Toeplitz operators in spaces of harmonic functions. Denote by
$\Hc(\O,\t)$ the subspace in $L_2(\O,\t)$, consisting of harmonic functions in a domain
$\O\subset\R^d$ and by $Q$ the orthogonal projection $Q: L_2(\O,\t)\to \Hc(\O,\t)$; this projection
is an integral operator with kernel $Q(x,y)$, $x,y\in\O$, the kernel being a harmonic function in
each variable $x$ and $y$. With a distribution $F$ having compact support in $\O$ we associate,
similarly to \eqref{distr.weight} the Toeplitz operator $T^\Hc_F:u\mapsto T^\Hc_Fu$,
$T^\Hc_Fu(z)=\langle F, Q(x,\cdot)u(\cdot)\rangle$. The expression for the action of the operator for
the case when $F$ is a Borel measure and the expressions for the sesquilinear form are analogous to
\eqref{meas}, \eqref{QForm}, \eqref{QFormMeas}. Similar to the case of analytical functions, we
associate with the distribution $F$ the matrix $\Hb_F$, with entries being $\langle F,f_\a
 \overline{f_\b}\rangle$, where $f_\a$ is some  system of harmonic polynomials in $\R^d$. Again, the rank
of the infinite matrix $\Hb_F$ does not exceed the rank of the operator $T^\Hc_F$. We, however, may
not include, as we have  done for analytical functions, the multiplicative functional parameter $g$,
since harmonic functions do not possess a multiplicative structure.

\section{Finite rank operators in dimension 1}\label{Sect3}
The aim of this section is to give a proof of the following result generalizing the Luecking theorem.
\begin{theorem}\label{Th:Dim1} Let $F$ be a distribution with compact support in the domain $\O\subset
\C^1.$ Suppose that the operator $T_F$ has finite rank $m$. Then there exist finitely many  points
$z_q\in\O$ $q=1,\dots,m_0$, $m_0\le m$, and differential operators $L_q=L_q(\partial x,\partial y),
q=1,\dots,m_0$ such that  $F=\sum L_q \d(z-z_q)$.\end{theorem}

We start with some observations about distributions in $\Ec'(\C)$. For such distribution we denote by
$\psupp F$ the complement of the unbounded component of the complement of $\supp F$.
\begin{lemma}\label{Lem.compsup}Let $F\in \Ec'(\C)$.
Then the following two statements are equivalent:\\
 a) there exists a distribution
$G\in \Ec'(\C)$ such that $\frac{\partial G}{\partial \bar{z}}=F$, moreover $\supp G\subset\psupp
F$;\\ b) $F$ is orthogonal to all polynomials of $z$ variable,  i.e. $\langle F,z^k\rangle=0$ for
all $k\in\Z_+$.
\end{lemma}
 \begin{proof} The implication $a)\Longrightarrow b)$ follows from the relation
\begin{equation}\label{byparts}
    \langle F, z^k\rangle=\langle \frac{\partial G}{\partial \bar{z}},z^k\rangle
    =\langle G ,\frac{\partial z^k}{\partial \bar{z}}\rangle=0.
\end{equation}

We prove that $b)\Longrightarrow a)$.  Put $G:=F*\frac{1}{\pi z}\in \Sc'(\C)$, the convolution being
well-defined since $F$ has compact support. Since $\frac{1}{\pi z}$ is the fundamental solution of
the Cauchy-Riemann operator $\frac{\partial}{\partial\bar{z}}$, we have $\frac{\partial G}{\partial
\bar{z}}=F$ (cf., for example, \cite{Ho}, Theorem 1.2.2).  By the ellipticity of the Cauchy-Riemann
operator, $\singsupp G\subset\singsupp F\subset \supp F$, in particular, this means that $G$ is a
smooth function outside $\psupp F$, moreover,  $G$ is analytic outside $\psupp F$ (by $\singsupp F$
we denote the singular support of the  distribution $F$, see, e.g., \cite{Ho}, the largest open set
where the distribution coincides with a smooth function).
  Additionally,
$G(z)=\langle F,\frac1{\pi(z-w)}\rangle =\pi^{-1}\sum_{k=0}^\infty z^{-k-1} \langle F,w^k\rangle=0$
if $|z|>R$ and $R$ is sufficiently large. By analyticity this implies $G(z)=0$ for all $z$  outside
$\psupp F$.\end{proof}

\begin{proof}[Proof of Theorem \ref{Th:Dim1}]

The  distribution in question $F$, as any distribution with compact support, is of finite order,
therefore it belongs to some Sobolev space, $F\in H^s$ for certain $s\in \R^1$. If $s\ge0,$ $F$ is a
function and must be zero by Luecking's theorem. So, suppose that $s<0$.

Consider the first $m+1$ columns in the matrix $\Ab_F$, i.e.
\begin{equation}\label{columns}
    a_{nk}=(T_Fz^k,z^l)=\langle F,z^k\bar{z}^l\rangle, l=0,\dots m;\  k=0,\dots.
\end{equation}
Since the rank of the matrix $\Ab_F$ is not greater than $m$, the columns are linearly dependent, in
other words, there exist coefficients $c_0,\dots,c_m$ such that $\sum_{l=0}^m a_{kl}c_l=0$ for any
$k\ge 0$. This relation can be written as
\begin{equation}\label{ColPolyn}
    \langle F,z^k h_1(\bar{z}) \rangle=\langle h_1(\bar{z})F, z^k\rangle=0, \ \ h_1(\bar{z}) = \sum_{k=0}^m c_l\bar{z}^l.
\end{equation}
Therefore the distribution $h_1(\bar{z})F\in H^s$ satisfies the conditions of Lemma \ref{Lem.compsup}
and hence there exists a compactly supported distribution $F^{(1)}$ such that $\frac{\partial
H^{(1)}}{\partial \bar{z}}=F_h$. By the ellipticity of the Cauchy-Riemann operator, the distribution
$F^{(1)}$ is less singular than $F$, $F^{(1)}\in H^{s+1}$. At the same time,
\begin{gather}\label{newmatrix}
    \langle F^{(1)}, z^k\bar{z}^l \rangle=
    (l+1)^{-1}\langle F^{(1)},\frac{\partial z^k\bar{z}^{l+1 }}{\partial \bar{z}}
    \rangle\\=(l+1)^{-1}\langle h(\bar{z})F, z^k\bar{z}^{l}\rangle =
    (l+1)^{-1}\langle F,
    z^k\bar{z}^{l}h(\bar{z})\rangle,
\end{gather}
and therefore the rank of the matrix $\Ab_{F^{(1)}}$ does not exceed the rank of the matrix $\Ab_F$.

We repeat this procedure sufficiently many (say, $ N=[-s]+1$) times and arrive at the distribution
$F^{(N)}$ in $L_2$, for which the corresponding matrix $\Ab_{F^{(N)}}$ has finite rank. By Luecking's
theorem, this may happen only if $F^{(N)}=0$.

Now we go back to the initial distribution $F$. Since, by our construction, $\frac{\partial
F^{(N)}}{\partial \bar{z}}=h_N(\bar{z})F^{(N-1)}$, we have that $h_N(\bar{z})F^{(N-1)}=0$ and
therefore $\supp F^{(N-1)}$ is a subset of the set of zeroes of the polynomial $h_N(\bar{z})$. On the
next step, since $\frac{\partial
F^{(N-1)}}{\partial \bar{z}}=h_{N-1}(\bar{z})F^{(N-2)}$,  we obtain that $\supp F^{(N-2)}$ lies in the union of  sets of zeroes of polynomials
$h_{N-1}(\bar{z})$ and $h_N(\bar{z})$. After having gone all the way back to $F$, we obtain that its
support is a finite set of points lying in the union of zero sets of polynomials $h_{j}$. A
distribution with such support must be a linear combination of $\d$ - distributions in these points
and their derivatives, $F=\sum L_q\d(z-z_q)$, where $L_q$ is some differential operator. Finally, to
show that the number of points $z_q$ does not exceed $m$, we construct for each of them the
interpolating polynomial $f_q(z)$ such that $L_q|f_q|^2\ne0$ at the point $z_q$   while at the points
$z_{q'},\ q'\ne q$, the polynomial $f_q$ has zero of sufficiently high order, higher than the order
of $L_{q'}$, so that $L_{q'}(f_q g)(z_{q'})=0$ for any smooth function $g$. With such choice of
polynomials, the matrix with entries $\langle F,f_q\overline{f_{q'}}\rangle $ is the diagonal matrix
with nonzero entries on the diagonal, and therefore its size (that equals the number of the points
$z_q$) cannot be greater than the rank of the whole matrix $\Ab_F$, i.e., cannot be greater than $m$.
\end{proof}

We note here that the attempt to extend the original proof of Luecking's theorem to the
distributional case would probably meet certain complications. Let us recall the crucial place in
\cite{Lue2}.

The matrix of the type \eqref{Matrix} is also considered, with a measure $\m$ standing on the place
of  the distribution $F$. Then, for a given $N$, the measure $\m^N=\otimes^N \m$ on $\C^N$ is
introduced, and  Lemma 5.1 is established, stating that if the Toeplitz operator $T_\m$ has rank
smaller than $N$, then for all symmetric polynomials $H_1(Z),H_2(Z)$ of the multi-dimensional complex
variable $Z=(z_1,z_2,\dots,z_N)\in \C^N$,
\begin{equation}\label{Luecking}
    \int H_1(Z)\overline{H_2(Z)}|V(Z)|^2d\m^N=0,
\end{equation}
where $V(Z)$ is the Vandermonde function, $V(Z)=\prod_{i<j}(z_i-z_j)$. To derive the finite rank
result from Lemma 5.1,  the following property is needed: the algebra generated by the functions of
the form $H_1(Z)\overline{H_2(Z)}$ is dense (in the sense of the uniform convergence on compacts) in the
space of symmetric continuous functions. This latter property is proved in \cite{Lue2} by an
ingenious reduction to the Stone-Weierstrass theorem.

Now, if $\mu=F$ is a distribution that is not a measure, the analogy of reasoning  in \cite{Lue2}
would require a similar density property, however not in the sense of the uniform convergence on
compacts, but in a stronger sense, the uniform convergence together with derivatives up to some fixed
order (depending on the order of the distribution $F$.)  The Stone-Weierstrass theorem seems  not to
help here since  it deals with uniform convergence only. Moreover, the required more general density
statement itself is \emph{wrong}, which follows   from the  construction below.

\begin{proposition}\label{example1} The algebra generated by the functions having  the form
$H_1(Z)\overline{H_2(Z)}$, where $H_1,H_2$ are symmetric polynomials of the variables
$Z=(z_1,\dots,z_N)$ is not dense in the sense  of the uniform  $C^l$-convergence  on compact sets  in
the space of $C^l$-differentiable symmetric functions, as long as $l\ge N(N-1)$.\end{proposition}
\begin{proof}We introduce the notations: $D_j=\frac{\partial}{\partial z_j}$, $\overline{D_j}=\frac{\partial}{\partial\bar{z_j}}$.
 Consider the differential operator $V(D)=\prod_{j<k}(D_j-D_k)$. It is easy to check that
$V(D)H$ is symmetric for any antisymmetric function $H(Z)$ and $V(D)H$ is antisymmetric for any
symmetric function $H(Z)$. Further on, consider any function $H(Z)$ of  the form
$H(Z)=H_1(Z)\overline{H_2(Z)}$ where $H_1(Z),H_2(Z)$ are analytic polynomials. If at least one of
them is symmetric, we have
\begin{equation}\label{sympol}
    V(D)V(\bar{D})H(0)=0.
\end{equation}
In fact,$ V(D)V(\bar{D}) H_1(Z)\overline{H_2(Z)} =[(V(D)H_1(Z)][\overline{V({D}){H_2(Z)}}]$. In the
last expression, for  the symmetric polynomial $H_l$ , the corresponding polynomial $V(D)H_l(Z)$ is
antisymmetric, and therefore equals zero for $Z=0$. Now consider the symmetric function
$|V(Z)|^2=V(Z)\overline{V(Z)}$. We have
 $$V(D)V(\bar{D})V(Z)\overline{V(Z)} =[V(D)V(Z)][V(\bar{D})\overline{V(Z)}].$$
Now note that $V(Z)=\sum_\kappa C_{\kappa}\prod z_j^{\kappa_j}$ where the summing goes over
multi-indices $\kappa=(\kappa_1,\dots,\kappa_N)$, $|\kappa|=N$ and not all of real coefficients
$C_\kappa$ are zeros. Simultaneously, $V(D)=\sum_{\kappa}C_{\kappa}\prod D_j^{\kappa_j}$ with the
same coefficients. We recall now  that $\prod D_j^{\kappa_j}\prod z_j^{\kappa'_j}=0$ if $\kappa\ne
\kappa'$  and it equals $\kappa!$ if $\kappa = \kappa'$. Therefore,
$V(D)V(Z)=\sum_{\kappa}C_{\kappa}^2\kappa!$ is a positive constant.  In this way we have constructed
the differential operator $V(D)V(\bar{D})$  of order $N(N-1)$, satisfying \eqref{sympol} for any
function of the form $H(Z)=H_1(Z)\overline{H_2(Z)}$ with symmetric $H_1,H_2$, and  not vanishing on
some symmetric differentiable function $|V(Z)|^2$. Therefore the function $|V(Z)|^2$ cannot be
approximated by linear combinations of the functions $H(Z)=H_1(Z)\overline{H_2(Z)}$ in the sense of
the uniform $C^{N(N-1)}$ convergence on compacts.\end{proof}

\section{The multi-dimensional case}\label{sect.Multi}

In this Section we extend our main Theorem \ref{Th:Dim1} to the case of Toeplitz operators in Bergman
spaces of analytical functions of several variables. For the case of a measure acting as weight,
there exist two proofs of this result, in \cite{BRChoe1} and \cite{RShir}. The first proof
generalizes the approach used in \cite{Lue2}, the other one uses the induction on dimension. As it
follows from Proposition \ref{example1}, for the case of distribution the approach of \cite{BRChoe1}
is likely to meet some complications. On the other hand, as we are going to show, the approach of
\cite{RShir} can be extended to the distributional case.
\begin{theorem}\label{Thm.Dim d} Let $F$ be a distribution in $\Ec'(\C^d)$. Consider the matrix
\begin{equation}\label{MatrtixDimd}
    \Ab_F=(a_{\a\b})_{\a,\b\in \Z_+^d}; \ a_{\a\b})=\langle F, Z^\a\bar{Z}^\b\rangle, \ Z=(z_1,\dots,z_d).
\end{equation}
Suppose that the matrix $\Ab_F$ has finite rank $m$. Then $\card\supp F\le m$ and $F=\sum L_q
\d(Z-Z_q)$, where $L_q$ are differential operators and $Z_q$ are some points in $\C^d$.
\end{theorem}
 We will perform the induction on dimension, proving a statement that is, actually, only formally
 weaker  than Theorem \ref{Thm.Dim d}, since, as
 it was explained in Sect.\ref{Sect2}, the rank of the matrix $\Ab_F$ does not grow if $F$
is replaced by $F_g$.
 \begin{proposition}\label{TheorDimDmore} Suppose that for any function $g(Z)$, analytic and bounded in a polydisk
 neighborhood of the
 support of the distribution $F$,  the conditions of Theorem  \ref{Thm.Dim d} are
 fulfilled with the distribution $F$ replaced  by $|g(Z)|^2F\equiv F_g$. Then $\card\supp F\le m$ and
  $F=\sum L_q
\d(Z-Z_q)$, where $L_q$ are differential operators.\end{proposition}

\begin{proof} For $d=1$ the statement of Proposition \ref{TheorDimDmore} coincides with the one of
Theorem \ref{Th:Dim1} that was proved in Sect.\ref{Sect3}. We suppose that we have established our
statement in dimension $d-1$ and consider the $d$-dimensional case. We denote the variables as
$Z=(z_1,Z'), \ Z'\in\C^{d-1}$.

 For a fixed function $g(Z)$ we denote by $G(g)=\pi_*F_g$ the distribution in $\Ec'(\C^{d-1})$ induced from $F_g$
  by the projection $\pi:Z\mapsto
 Z'$: for $u\in C^\infty(\C^{d-1})$
 \begin{equation}\label{proj}
    \langle G(g), u\rangle = \langle F_g, 1_{\C^1}\otimes u\rangle.
 \end{equation}
Although the function $g$ is defined only in a polydisk, the expression in \eqref{proj} is well
defined since this polydisk contains $\supp F$.

Consider the submatrix $\Ab'_{F_g}$ in the matrix $\Ab_{F_g}$ consisting only of those
$a_{\a\b}=\langle |g|^2 F, Z^a\overline{Z}^\b\rangle$  for which $\a_1=\b_1=0$. It follows from
\eqref{proj}, that the matrix $\Ab'_{F_g}$ coincides with the matrix $\Ab_{G(g)}$ constructed for the
distribution $G(g)$ in dimension $d-1$. Thus, the matrix $\Ab_{G(g)}$, being a submatrix of a finite
rank matrix, has a finite rank itself, moreover, $\rank\Ab_{G(g)}\le m$. By the inductive assumption,
this implies that the distribution $G(g)$ has  finite support consisting of $m(g)\le m$ points
$\z_1(g),\dots,\z_{m(g)}$;  $\z_q(g)\in \C^{d-1}$ (the notation reflects the fact that both the points
and their quantity may depend on the function $g$). Among all functions $g$, we can find the one,
$g=g_0$, for which $m(g)$ attains its maximum value $m_0\le m$. Without losing in generality, we can
assume that $g_0=1$.

Fix an $\e>0$, sufficiently small, so that $2\e$-neighborhoods of $\z_q(1)$ are disjoint, and
consider the functions $\varphi_q(z')\in C^{\infty}(\C^{d-1})$, $q=1,\dots, $ such that $\supp
\varphi_q $ lies in the $\e$-neighborhood of the point $\z_q(1)$ and $\varphi(z')=1$ in the
$\frac\e2$-neighborhood of $\z_q(1)$.  We fix an analytic function $g(z)$ and consider for any $q$
the distribution $\F_q(t,g)\in \Ec'(\C^d),$ $\F_q(t,g)=|1+tg|^2 \varphi_q(Z')F=\varphi_q(Z')F_{1+tg}$.
For $t=0$, $\F_q(t,g)=\varphi_q(Z')F$,  the point $\z_q(1)$ belongs to the support of $\pi_*
\F_q(0,g)$, and therefore for some function $u\in C^\infty(\C^{d-1})$, $\langle\pi_*
\F_q(0,g),u\rangle\ne0$. By continuity, for $|t|$ small enough, we still have $\langle\pi_*
\F_q(t,g),u\rangle\ne0$, which means that the $\e$-neighborhood of the point $\z_q(1)$ contains at
least one point in the support of the distribution $G(1+tg)$.  Altogether, we have not less than
$m_0$ points of the support of $G(1+tg)$ in the union of $\e$-neighborhoods of the points $\z_j(1)$.
However, recall, the support of $G(1+tg)$ can never contain more than $m_0$ points, so we deduce that
for $t$ small enough, there are no points of the support of $G(1+tg)$ outside the $\e$-neighborhoods
of the points $\z_q(1)$, so
\begin{equation}\label{nbhood}
\supp G(1+tg)\cap\{Z':|Z'-\z_q|>\e\}=\varnothing\end{equation} for $|t|$ small enough (depending on
$g$.)
 Now we introduce a function $\p\in C^\infty(\C^{d-1})$ that  equals $1$ outside
$2\e$-neighborhoods of the points $\z_q(1)$ and vanishes in $\e$-neighborhoods of these points. By
\eqref{nbhood}, the distribution $\p G(1+tg)$ equals zero for any $g$, for $t$ small enough. In
particular, applying this distribution to the function $u=1$, we obtain
\begin{equation}\label{Ftg}
\langle\p G(1+tg),1\rangle=\langle\p F,|1+tg|^2\rangle=\langle \p F, 1+2t\re g+ t^2|g|^2\rangle=0.
\end{equation}
By the arbitrariness of $t$ in a small interval, \eqref{Ftg} implies that $\langle \p F,
|g|^2\rangle=0$ for any $g$. Now we take $g$ in the form $g=g_1+g_2$, where $g_1,g_2$ are again
functions analytical in a polydisk neighborhood of $\supp F$. Then we have
    $$\langle\p F, |g_1|^2+2\re(g_1\overline{g_2})+|g_2|^2\rangle=\langle\p F, 2\re(g_1\overline{g_2})
    \rangle=0.$$
Replacing here $g_1$ by $ig_1$, we obtain $\langle\p F, 2\im(g_1\overline{g_2})
    \rangle=0$, and thus
    \begin{equation}\label{Fgg}
\langle\p F, g_1\overline{g_2}
    \rangle=0.
    \end{equation}

    Any  polynomial $p(Z,\bar{Z})$ can be represented as a linear combination of functions
    of the form $g_1\overline{g_2}$, so, \eqref{Fgg} gives
    \begin{equation}\label{FP}
\langle\p F, p(Z,\bar{Z})
    \rangle=0.
    \end{equation}
Now  we take any function $f\in C^\infty(\C^d)$ supported in the neighborhood $V$ of $\supp F$ such
that $f=0$  on the support of $\p$. We can approximate $f$
 by polynomials of the form $p(Z,\bar{Z})$ uniformly on $\overline{V}$ in the sense of $C^l$, where
 $l$ is the order of the distribution $F$. Passing to the limit in \eqref{FP}, we obtain $\langle\p F,
 f \rangle=\langle F,
 f \rangle=0.$

 The latter relation shows that $\supp F\subset \bigcup_q\{Z: |Z'-\z_q(1)|<2\e\}$. Since $\e>0$ is arbitrary,
 this  implies that  $\supp F$
lies in the union of affine subspaces $Z'=\z_j$, $j=1,\dots,m_0$ of complex dimension $1$. Now we
repeat the same reasoning having chosen instead of $Z=(z_1,Z')$ another decomposition of the complex
variable $Z$: $Z=(Z'',z_d)$. We obtain that for some points $\x_k\in \C^{d-1}$, no more than $m$ of
them, the support of $F$ lies in the union of subspaces $Z''=\x_k$. Taken together, this means that,
actually, $\supp F$ lies in the intersection of these two systems of subspaces, which consists of no
more than $m^2$ points $Z_s$. The number of points is finally reduced to $m_0\le m$ in the same way
as in Theorem\ref{Th:Dim1}, by choosing a special system of interpolation functions.
\end{proof}

\section{Harmonic functions}\label{Sect:Harm}
The aim of this section is to establish finite rank results for  Toeplitz operators corresponding to
the Bergman  spaces of harmonic functions. The main difference with the analytical case lies in the
circumstance that the space of harmonic functions does not possess the multiplicative  structure.
Therefore, in the process of dimension reduction, similar to the one we used in the proof of Theorem
\ref{Thm.Dim d}, we are not able to introduce the functional parameter (denoted by $g$ there.) As a
result of this circumstance, we can prove the finite rank theorem only in the case of $F$ being a
measure and not a more singular distribution. In order to justify this shortcoming, we conclude the
section by presenting an example of a  singular distribution  with rather large support (and thus
non-discrete), that projects to a discrete measure, whatever the direction of the projection. Thus, a
considerable part of $F$ becomes invisible after being projected. This example, although not
contradicting directly the finite rank property, indicates that the reduction of dimension might be
not sufficient to prove the result.

We start with the even-dimensional case. Here the problem with harmonic spaces reduces easily to the
analytical case (in fact, we could have used a reference to \cite{BRChoe1} instead).
\begin{theorem}\label{ThHarmEven} Let $d=2n$ be an even integer. Suppose that for a
certain distribution $F\in \Ec'(\R^n)$ the matrix $\Hb_F$ defined in Section \ref{Sect2} has rank
$m<\infty$. Then the distribution $F$ is a sum of $m_0\le m$ terms, each supported in one point:
$F=\sum L_j\d(x-x_q)$, $x_q\in \R^n$, $L_q$ are differential operators in $\R^n$.\end{theorem}
\begin{proof} We identify the space $\R^d$ with the complex space $\C^n$. Since the functions
$z^\a,\bar{z}^\b$ are harmonic, the matrix $\Ab_F$ can be considered as a submatrix of $\Hb_F$, and
therefore it has rank not greater than $m.$ It remains to apply Proposition \ref{TheorDimDmore} to
establish that the distribution $F$ has the required form, with no more  than $m$ points $x_q$.
\end{proof}

The odd-dimensional case requires considerably more work. We will use again  a kind of dimension
reduction, however, unlike the analytic case,  we will need projections of the distribution to
one-dimensional subspaces.

 Let $S$ denote the unit sphere in $\R^d$, $S=\{\z\in\R^d:|\z|=1\}$ and let $\s$ be the Lebesgue
measure on $S$. For $\z\in S$, we denote by $\Lc_\z$ the one-dimensional subspace in $\R^d$ passing
through $\z$, $\Lc_\z=\z \R^1$. For a distribution $F\in \Ec'(\R^d)$, we define the distribution
$F_\z\in\Ec'(\R^1)$ by setting $\langle F_\z,\f\rangle=\langle F,\f_z\rangle$, where $\f_z\in
C^\infty(\R^d)$ is $\f_z(x)=\f(x\cdot z)$. The distribution $F_\z$ can be understood as result of
projecting of $F$ to $\Lc_\z$ with further transplantation of the projection, $\pi_*^{\Lc_\z}F$, from
the line $\Lc_\z$ to the standard line $\R^1$.  The  Fourier transform $\Fc F_\z$ of $F_\z$ is
closely related with $\Fc F$:
\begin{equation}\label{ProjDistr}
    \Fc (F_\z)(t)=(\Fc
F)(t\z).
\end{equation}

Further on, we will restrict ourselves to the case when the distribution $F$ is a finite complex
Borel measure  $\m$. Here we will use the notation $\m_\z$ instead of $F_\z$ .

 We need to recall  certain facts in harmonic analysis. In the one-dimensional
case, they were  proved by N.Wiener as long ago as in 1919; the multi-dimensional version seems to be
folklore, however the formulations we found in the literature, see \cite{mattila}, are slightly
weaker than the ones we need.

Let $\m$ be a finite complex  Borel measure in $\R^d$. We define
 $$\lfloor\mu\rfloor=\left(\sum_{\x\in\R^d}|\m(\{\x\})|^2\right)^{\frac12}.$$
 Of course, $\lfloor\mu\rfloor$ is finite for a finite measure and it vanishes if and
  only if $\m$ has no atoms.

 \begin{lemma}\label{lemWiener}Let $\m$ be a finite Borel measure $\R^d$ and $h$ be a function in
 $L_1(\R^1)$. Denote by $\Fc\m$ the Fourier transform of $\m$. Then
 \begin{equation}\label{Wiener1}
    \lim_{R\to\infty}R^{-d}\int_{\R^d}h(R^{-1}\x)\Fc\m(\x)d\x=\m(\{0\})\int_{\R^d}h(\x)d\x.
 \end{equation}
\end{lemma}
\begin{proof}By Plancherel identity,
we have
$$\lim_{R\to\infty}R^{-d}\int_{\R^d}h(R^{-1}\x)\Fc\m(\x)d\x=\lim_{R\to\infty}\int_{\R^d}(\Fc h)(R x)d\m(x).$$
Now note that $\Fc h(0)=\int_{\R^d}h(\x)d\x$ and $\lim_{\R\to0}(\Fc h)(Rx)=0$ for $x\ne0$ by
Riemann-Lebesgue lemma. The proof completes by applying the Lebesgue dominant convergence theorem.
\end{proof}
\begin{corollary}\label{CorWiener} Under the conditions of Lemma \ref{lemWiener},
\begin{equation}\label{Wiener2}
\lim_{\R\to\infty}R^{-d}\int_{\R^d}h(R^{-1}\x)|\Fc\m(\x)|^2
d\x=\lfloor\m\rfloor^2\int_{\R^d}h(\x)d\x.\end{equation}
\end{corollary}
\begin{proof}We define the measure $\check{\m}$ as $\check{\m}(E)=\overline{\m(-E)}$ for any Borel
set $E$ and introduce $\n=\m*\check{\m}$. Then $\Fc\n=|\Fc\m|^2$ and $\n(0)=\lfloor\m\rfloor^2$. It
remains to apply Lemma \ref{lemWiener} to the measure $\m$.
\end{proof}

We are going to use Corollary \ref{CorWiener} to relate the properties of the family of measures
$\m_\z,\ \z\in S$, with the properties of $\m$.
\begin{lemma}\label{LemCont} Let $\m$ be a finite compactly supported complex Borel measure on
$\R^d$. Then the following two statements are equivalent:\\
a) The measure $\m$ is continuous, i.e., $\m(\{x\})=0$ for any $x\in \R^d$,\\
 The measure $\m_\z$ is continuous for $\s$-almost all $\z\in S$.
\end{lemma}
\begin{proof}
We take a function $h(\x)$, depending only on $|\x|$, $h(\x)=H(|\x|)$ such that
$\int_{\R^d}h(\x)d\x=1$. So, $\int_{\R} |r|^{d-1}H(|r|)dr=\frac{2}{\s(S)} $. By Corollary
\ref{CorWiener}, used in dimension $1$ for $\m_\z$,
\begin{equation*}
    \lfloor\m_\z\rfloor^2=\lim_{R\to\infty}\frac{\s(S)}{2R}\int_{\R}|R^{-1}r|^{d-1}H(R^{-1}|r|)|(\Fc\m_\z)(r)|^2dr.
\end{equation*}
In what follows we apply the Lebesgue dominant compactness theorem to justify the passing to a limit:
\begin{gather}\label{WienerLONG}\nonumber
    \frac{1}{\s(S)}\int_S\lfloor\m_\z\rfloor^2 d\s(\z)\\
    =\int_S
    \lim_{R\to\infty}\frac{1}{2R}\int_{\R}|R^{-1}r|^{d-1}H(R^{-1}|r|)|(\Fc\m_\z)(r)|^2drd\s(\z)\nonumber\\
    =\lim_{R\to\infty}\frac{1}{2R^d}\int_S\int_{\R}|r|^{d-1}H(R^{-1}|r|)|(\Fc\m_\z)(r)|^2drd\s(\z)\nonumber\\
=\lim_{R\to\infty}\frac{1}{R^d}\int_S\int_0^\infty r^{d-1}H(R^{-1}r)|(\Fc\m)(r\z)|^2drd\s(\z)\nonumber\\
=\lim_{R\to\infty}\frac{1}{R^d}\int_{\R^d}h(R^{-1}\x)|(\Fc\m)(\x)|^2d\x=\lfloor\m\rfloor^2.
\end{gather}
Hence, $\lfloor\m\rfloor=0$ if and only if $\lfloor\m_\z\rfloor=0$ for almost all $\z\in S$.
\end{proof}
\begin{corollary}\label{WienerLongCor} For a finite complex Borel measure $\m$ with compact support in $\R^d$
the following three statements are equivalent:\\
a) $\m$ is discrete;\\
 b) $\m_\z$ is discrete for all $\z\in S$;\\
 c) $\m_\z$ is discrete for $\s$-almost all $\z\in S$.
\end{corollary}
\begin{proof} The implications $a)\Longrightarrow b)$ and $b)\Longrightarrow c)$ are obvious. To
establish $c)\Longrightarrow a)$, we denote by $\m^c$ the continuous part of $\m$. Then the statement
$c)$ means that $(\m^c)_\z$ is discrete for $\s$-almost all $\z\in S$. On the other hand, by Lemma
\ref{LemCont} applied to $\m^c$, the measure $(\m^c)_\z$ is continuous for $\s$-almost all $\z\in S$.
Being both discrete and continuous, the measure $(\m^c)_\z$ is zero for  $\s$-almost all $\z\in S$.
Passing to the Fourier transform, we obtain $(\Fc\m^c)(r\z)=0$ for all $r$ for $\s$-almost all $\z\in
S$. Now, since the Fourier transform $\Fc\m^c$ is smooth, this means that $\m^c=0.$
\end{proof}

Now we return to our finite rank problem.
\begin{theorem}\label{TheoHarm}
Let $d\ge3$ be an odd integer, $d=2n+1$. Let $\m$ be a finite complex Borel measure in $\R^d$ with
compact support. Suppose that the matrix $\Hb_\m$ has finite rank $m$. Then $\supp \m$ consists of no
more than $m$ points.
\end{theorem}

\begin{proof}Fix some $\z\in S$ and chose some $d-1=2n$-dimensional linear subspace $\Lc\subset \R^{d}$ containing $\Lc_\z$.
 We choose  the co-ordinate system
$x=(x_1,\dots,x_d)$ in $\R^d$ so that the subspace $\Lc$ coincides with $\{x:x_d=0\}$. The
even-dimensional real space $\Lc$ can be considered as the $n$-dimensional complex space $\C^n$ with
co-ordinates $z=(z_1,\dots,z_n)$,  $z_j=x_{2j-1}+ix_{2j}$, $j=1,\dots,n$. The functions
$(z,x_d)\mapsto z^\a$, $(z,x_d)\mapsto \bar{z}^\b$, $\a,\b\in (\Z_+)^d$, are harmonic polynomials in
$\C^d\times  \R^1$. Moreover, by definition,
$\langle\m,z^\a\bar{z}^\b\rangle=\langle\pi^{\C^n}_*\m,z^\a\bar{z}^\b \rangle$. Hence,  the matrix
$\Ab_{\pi^{\C^n}_*\m}$ is a submatrix of the matrix $\Hb_{\m}$, and the former has not greater rank
than the latter, $\rank(\Ab_{\pi^{\C^n}_*\m})\le m$. So we can apply Theorem \ref{Thm.Dim d} and
obtain that the measure $\pi^{\C^n}_*\m$ is discrete and its support contains not more than $m$
points. Now we project the measure $\pi^{\C^n}_*\m$ to the real one-dimensional  linear subspace
$\Lc_\z$ in $L$.  We obtain the same measure as if we had projected $\m$ to $\Lc_\z$ from the very
beginning, and not in two steps i.e.,  $\pi^{\Lc_\z}_*\m$. As a projection of a discrete measure,
$\pi^{\Lc_\z}_*\m$ is discrete and has no more than $m$ points in the support. By our definition of
the measure $\m_\z$ as $\pi^{\Lc_\z}_*\m$ transplanted to $\R^1$, this means  that $\m_\z$ is
discrete.

Due to the arbitrariness of the choice of $\z\in S$, we obtain that all measures $\m_\z$ are
discrete. Now we can apply Corollary \ref{WienerLongCor} and obtain that the measure $\m$ is discrete
itself. Finally, in order to show that the number of points in $\supp \m$ does not exceed $m$, we
chose $\z\in S$ such that no two points in $\supp \mu$ project to the same point in $\Lc_\z$. Then
the point masses of $\m$ cannot cancel each other under the projection, and thus $\card\supp
\m=\card\supp\m_\z\le m$.

The number of points in the support of $\m$ is estimated in the same way as in Theorems \ref{Th:Dim1}
and \ref{Thm.Dim d}.
\end{proof}
 The analysis of the reasoning in the proof shows that the only essential obstacle for extending Theorem
\ref{TheoHarm} to the case of distributions is the limitation set by Corollary \ref{WienerLongCor}.
If we were able to prove this Corollary for distributions, all other steps in the proof of Theorem
\ref{TheoHarm} would go through without essential changes.  However, it turns out that not only the
proof of Corollary \ref{WienerLongCor} cannot be carried over to the distributional case, but,
moreover, the Corollary itself becomes wrong. The example that we present does not disprove Theorem
\ref{TheoHarm} for distributions, however it indicates that the proof, if exists, should involve some
other ideas.

\begin{example}\label{example}
Let $d\ge2$. We consider the Schwartz distribution $F\in \Sc(\R^d)$ that has $\cos|\x|$ as its
Fourier transform. By the Paley-Wiener theorem, since $\Fc F$ is an entire function of exponential
type, $F$ has compact support, $F\in\Ec'(\R^d)$. By \eqref{ProjDistr} and spherical symmetry, for any
$\z\in S$, $F_\z=\Fc^{-1}(\cos \t)=\frac12 (\d_1+\d_{-1})$. If $F$ were a measure, then, by Corollary
\ref{WienerLongCor} it would be discrete.  This, however is impossible since $F$, together with $\Fc
F$, is rotationally invariant; being both discrete  and rotationally invariant, $F$ must have support
in the origin, which contradicts the above expression for $F_\z$. The construction also shows that $F$
is the unique distribution that has $\frac12 (\d_1+\d_{-1})$ as its one-dimensional projections.  Of
course, we could have  directly checked that $F$ is not a measure, using the fact that $F$ is,
actually, the solution $u(x,t)$, $t=1$, for the wave equation $u_{tt}-\D_x u=0$ with initial
conditions $u(\cdot,0)=\d$, $u_t(\cdot,0)=0$. Moreover, from the classical Poisson formulas it
follows that
  $\supp F$ is the sphere $\{|x|=1\}$ for odd $d$ and the ball  $\{|x|\le1\}$  for even $d$. Note,
  however, that in neither dimension $F$ generates a finite rank Toeplitz operator.
\end{example}
\section{Discussion}\label{discussion}
In the process of exploring the finite rank conjecture, a number of interesting open questions
 arise. The case of analytical functions is studied completely. However, in the case of harmonic
 functions the finite rank conjecture is open for weights being distributions that are not measures.
The complete solution of this problem would follow from the positive answer to the next  question.
Let $d\ge3$, $F\in\mathcal{E}'(\R^d)$. Suppose that $\pi_*^H F$ is a distribution with a finite
support for every subspace $H\subset\R^d$ with $\dim H=d-1$. Is it true that the support of $\mu$ is
finite? As Example~\ref{example} shows, the answer is negative, if we consider subspaces of dimension
$1$  instead.

Further possible versions of the finite rank conjecture may involve some other elliptic equations
playing the part of the Cauchy-Riemann or the Laplace equations in the problem. The first interesting
candidate for the study here is the Helmholtz  operator $H_E u=\D u+Eu, \ E>0$. Let $P_{\Hc_E}$ be
the orthogonal projection from $L_2(\O)$ to the subspace $\Hc_E(\O)$ consisting of solution of the
Helmholtz equation. With a function (or a compactly supported  distribution) $F$ we associate the
Toeplitz operator $T_F: u\mapsto P_{\Hc_E}uF , \ u\in \Hc_E(\O)$. Which restrictions on $F$ are
imposed by the the condition that the  operator $T_F$ has a finite rank?  The question is of  a
certain importance for the scattering theory. It is easy to show that if $T_F$ is zero then $F$ must
be zero. However it is unclear at the moment how to handle the case of a positive rank. For the
Toeplitz operator corresponding to the projection onto the subspace of solutions of a general
elliptic equation, even the case of rank $0$ is unresolved.


\begin{thebibliography}{333}

\bibitem{BRChoe1} Choe, B.R., On higher dimensional Luecking's theorem.
\verb"http://math.korea.ac.kr/~choebr/papers/luecking.pdf"


\bibitem{BRChoe2}Choe, B.R.,  Koo, H.,  Young J.L.
 Finite sums of Toeplitz products on the polydisk,\\
 \verb"http://math.korea.ac.kr/~choebr/papers/finitesum_polydisk.pdf"




\bibitem{Cuc}\u{C}u\u{c}kovic, Z., Louchihi, I.,
Finite rank commutators and semicommutators of Toeplitz operators with quasihomogeneous symbols.
\verb"http://www.math.utoledo.edu/~ilouhic2/PDF/rank.pdf"






\bibitem{Vas}Grudsky, S. M.; Vasilevski, N. L. Toeplitz operators
on the Fock space: radial component effects.
    Integral Equations Operator Theory,  44  (2002),  no. 1, 10--37.
\bibitem{Ho} H\"ormander, L. ; The analysis of linear partial differential operators. V.1.Springer.,
1983.
\bibitem{Le} Le, T. A refined Luecking's theorem
 and finite-rank products of Toeplitz operators, arXiv:0802.3925
 \bibitem{Le1} Le T. Finite rank products of Toeplitz operators in several complex
 variables. {\small\verb"http://individual.utoronto.ca/trieule/Data/FiniteRankToeplitzProducts.pdf"}.
\bibitem{Lue2}Luecking, D.
Finite rank Toeplitz operators on the Bergman space,
Proc. Amer. Math. Soc. 136 (2008), no. 5, 1717--1723.
\bibitem{mattila} Mattila, P. Spherical averages of Fourier transforms of measures with finite energy;
 dimension of intersections and distance sets,  Mathematika 34(2) (1987), 207--228.

\bibitem{RaiWar}Raikov, G. Warzel, S.,
Quasi-classical versus non-classical spectral
    asymptotics for magnetic Schr\"odinger operators
     with decreasing electric potentials.  Rev. Math. Phys.
      14  (2002),  1051--1072.

      \bibitem{RShir} Rozenblum, G., Shirokov, N., Finite rank Bergman-Toeplitz
       and Bargmann-Toeplitz operators in many dimensions, arXiv:0802.0192


\end{thebibliography}
\end{document}